\documentclass[12pt]{article}
\def\dateref{09 April 2019}
\usepackage{amssymb, amsmath, amsthm, amscd, graphicx, dsfont, mathabx}
\usepackage[T1]{fontenc}
\usepackage[all]{xy}
\usepackage{pictexwd,dcpic}
\newtheorem{lem}{Lemma}[section]
\newtheorem{pro}{Proposition}[section]
\newtheorem{thm}{Theorem}[section]

\newtheorem{defn}{Definition}[section]

\def\be{\begin{equation}}
\def\ee{\end{equation}}

\def\be{\begin{equation}}
\def\ee{\end{equation}}

\def\cE{{\mathcal E}}

\def\cL{{\mathcal L}}

\def\N{\mathbb{N}}
\def\Z{\mathbb{Z}}
\def\R{\mathbb{R}}

\def\B{\mathbb{B}}

\def\tA{{\tt A}}

\def\tb{{\tt b}}

\def\brakt#1#2{\langle\,#1,#2\,\rangle}

\def\colon{\,{:}\;}
\def\prf{\medbreak\noindent{\bf Proof.}\enspace}
\def\qed{\hspace*{\fill}\hbox{\vrule height 7pt \kern-.3pt
     \vbox{\hrule width 7pt
     \kern6.6pt\hrule width 7pt }\kern-.3pt\vrule height 7pt
     }\par}
\def\ra{\rightarrow}
\def\cirk{{\scriptstyle\circ}}
\def\!{\mskip-\thinmuskip}

\def\!{\mskip-\thinmuskip}

\begin{document}

\iffalse
\\title{\bf On the equivalence of weak Gibbs measures and equilibrium measures
for shift spaces}
\fi
\title{\bf Weak Gibbs and Equilibrium Measures  for Shift Spaces}

\author{C.-E. Pfister\footnote{E-mail: charles.pfister@epfl.ch}\\
Faculty of Basic Sciences,
Section of Mathematics,\\ EPFL,
Station 8, CH-1015 Lausanne, Switzerland
\and W.G. Sullivan\footnote{E-mail:
Wayne.Sullivan@ucd.ie}\\
               Department of Mathematics and Statistics  UCD,\\
               Belfield, Dublin 4, Ireland}

 \date{\dateref}

\maketitle

\def\tb{{\tt b}}
\noindent
{\bf Abstract\colon}
For a large class of irreducible shift spaces $X\subset\tA^{\Z^d}$, with $\tA$  a finite alphabet,
and for absolutely summable potentials $\Phi$, we prove that  equilibrium measures for $\Phi$ are weak Gibbs measures.
In particular, for $d=1$, the result holds for irreducible sofic shifts.
%\newpage

%\vspace*{1.5cm}
\section{Introduction}

Equilibrium measures and Gibbs measures are  basic concepts which occur
 naturally in statistical mechanics and dynamical systems. For shift spaces $X$ (or lattice models) with a finite alphabet $\tA$ and for absolutely convergent potentials $\Phi$ Ruelle proved that Gibbs measures and equilibrium measures are equivalent provided that $X\subset \tA^{\Z^d}$ is
a subshift of finite  type \emph{and} verifies a supplementary condition (condition $(D)$, \cite{Ru}  chapter 4) see
comment after definition \ref{defnperiodicdecoupling}. In this note
we consider shift spaces verifying only a condition weaker than condition $(D)$ of Ruelle, called below \emph{decoupling condition},
and  prove that equilibrium measures are equivalent to weak Gibbs measures (in the sense of dynamical systems theory).
When $d=1$ we can weaken the decoupling condition and obtain the same result for a larger class of shift spaces. Our results apply for example to all irreducible sofic shifts.
The equivalence of equilibrium measures with weak Gibbs measures implies that
the empirical measures on the probability space $(X,\mu)$, with $\mu$ an equilibrium measure, have good large deviations estimates \cite{PS}. As a consequence of the contraction principle one has good large deviations estimates for macroscopic observables, which are defined as  averages of local functions with respect to the
shift action.

In the next section the precise setting and the results are formulated. The main theorem is proved
in the last section.

\section{Setting and main result}

Let $L=\Z^d$ and $\tA$ be a finite set (with discrete topology). The full shift
$\tA^L$ is a compact (metric) space for the product topology. The natural action by translation of $L$ on
 $\tA^L$  is denoted by $T\colon L\times \tA^L\ra \tA^L$,
$$
(j,x)\mapsto T^jx\quad\text{where}\quad (T^jx)(k):=x(k+j)\,,
$$
$x(k)$ being the $k^{th}$ coordinate of $x$. A shift space $X$ is a closed $T$-invariant subset of $\tA^L$.
$C(X)$ is the set of continuous functions on $X$ with the sup-norm $\|\cdot\|_\infty$,
$M_1(X)$ is  the set of Borel probability measures on $X$ (with
the topology of  weak convergence) and
$M_1(X,T)$ the subset of $T$-invariant probability measures.

For $\Lambda\subset L$ we set $\Lambda^c:=L\backslash\Lambda$ and $|\Lambda|$ denotes the cardinality of $\Lambda$ when $\Lambda$ is finite. For each integer $m$
$$
\Lambda_m:=\{i\in  L\colon \max\{|i_k|\colon k=1,\ldots,d\}\leq m\}\,.
$$
For $d=1$, we use also the notations
$$
[m_1,m_2]=\{i\in\Z\colon m_1\leq i\leq m_2]\,,\, (-\infty,m)=\{i\in\Z\colon i< m\}\,,\,
(m,\infty)=\{i\in\Z\colon i> m\}\,.
$$
We use $J_\Lambda$ for the projection map
$$
J_\Lambda\colon \tA^{\Lambda^\prime}\ra \tA^\Lambda\,,\;
J_\Lambda(x):=(x(k)\colon k\in\Lambda)\quad \text{for $\Lambda\subset\Lambda^\prime\subset L$}\,.
$$
Let $\Lambda\subset L$, $|\Lambda|<\infty$.
The empirical measure $\cE_\Lambda(x)$  is the discrete measure
$$
\cE_\Lambda(x):=\frac{1}{|\Lambda|}\sum_{j\in\Lambda}\delta_{T^jx}\quad\text{and}\quad
\brakt{\cE_\Lambda(x)}{f}:=\int f\,d\cE_\Lambda(x)\,.
$$
We also write  $J_m$ for $J_{\Lambda_m}$,
$X_n$ for $J_n(X)$ and
$\cE_m(x)$ for $\cE_{\Lambda_m}(x)$.

An absolutely summable potential $\Phi=\{\Phi_A\}$ is a family of continuous functions $\Phi_A:X\ra\R$ indexed by the finite subsets of $L$, such that $\Phi_A$ is
a local function, that is $\Phi_A(x)=\Phi_A(y)$ whenever
$J_A(x)=J_A(y)$, $\Phi_A\cirk T^a=\Phi_{A+a}$ and $\|\Phi\|:=\sum_{A\ni 0}\|\Phi_A\|_\infty<\infty$.
The set of absolutely summable potentials $\Phi$ with the norm $\|\Phi\|$ is a Banach space $\B$.
To each $\Phi$ we associate a continuous function
$$
\varphi_\Phi:=\sum_{A\ni 0}\frac{\Phi_A}{|A|}\,.
$$
Let $\Lambda$ and $M$ be disjoint subsets of $L$, $\Lambda$ finite. We set
$$
U_\Lambda:=\sum_{A\subset\Lambda}\Phi_A\quad\text{and}\quad W_{\Lambda,M}:=\!\!\!\!\!\!\!\!
\sum_{\substack{A\subset \Lambda\cup M\colon\\A\cap\Lambda\not=\emptyset\,A\cap M\not=\emptyset}}
\!\!\!\!\!\!\Phi_A\,.
$$
The pressure for a continuous function $\psi$ is denoted by $P(\psi)$ (\cite{Ru}, \cite{Wa}).
If $\psi=\varphi_\Phi$, with $\Phi\in\B$, then  $P(\varphi_\Phi)=\lim_{n\ra\infty}P_n(\Phi)$, where
$$
P_n(\Phi):=\frac{1}{|\Lambda_n|}\ln\sum_{x\in X_n}\exp U_{\Lambda_n}(x)\,.
$$
Let
$$
B_m(x):=\{y\in X\colon J_{\Lambda_m}(y)=J_{\Lambda_m}(x)\}\,.
$$

\begin{defn}\label{defnweakgibbs}
A probability measure $\nu\in M_1(X)$ is  a \emph{weak Gibbs measure for the continuous function $\psi$} if for any $\delta>0$ there exists $N_\delta$ such that
\begin{equation}\label{esti}
\sup_{x\in X}\Big|\frac{1}{|\Lambda_m|}\ln\nu(B_{m}(x))-\brakt{\cE_m(x)}{\psi}\Big|\leq \delta\quad\forall\,m\geq N_\delta\,.
\end{equation}
\end{defn}

The set of weak Gibbs measures for $\psi$ is convex. Indeed, if $\nu_1$ and $\nu_2$ verify \eqref{esti},
then this is also true for $\min_{i=1,2}\nu_i(B_{m}(x))$ and $\max_{i=1,2}\nu_i(B_{m}(x))$ in place of
$\nu(B_{m}(x))$. Hence \eqref{esti} holds for $a\nu_1+(1-a)\nu_2$, $0<a<1$.
If $\nu\in M_1(X,T)$ is a weak Gibbs measure for $\psi$, then $P(\psi)=0$ and $\nu$ is an
equilibrium measure for $\psi$ (see \cite{PS}).

\noindent
{\bf Remark.\,}
Definition \ref{defnweakgibbs} is equivalent to that of
\cite{PS} when $X$ is a shift space.
The terminology is that in usage in dynamical systems theory. In statistical physics
there is another notion  of ``weak Gibbs measure''. See e.g. \cite{Le}.
\qed

\begin{defn}\label{defntangentfunctional}
A probability measure $\nu\in M_1(X,T)$ is a \emph{tangent functional to the pressure $P$ at the continuous function $\varphi$} if it is an element of
$$
\partial P(\varphi):=\Big\{\nu\in M_1(X,T)\colon P(\varphi+f)\geq P(\varphi)+\int f\,d\nu\,,\;\forall f\in C(X)\Big\}\,.
$$
\end{defn}

In the present setting equilibrium measures and  tangent functionals to the pressure are equivalent  (see \cite{Wa}  theorems 9.15 and  8.2).

\begin{defn}\label{defnperiodicdecoupling}
A shift space $X\subset \tA^L$ satisfies the \emph{decoupling condition} if \\
1) there exists a function $q\colon \N\to\N$, written $q_m$,
such that $\lim_{m\to\infty} q_m/m=0$;\\
2) for $m\in \N$ and  $x,y\in X$ there
exist $z\in X$ and $\ell\in\Lambda_{q_m}$ such that
$$
J_{\Lambda_m}(T^{-\ell}z)=J_{\Lambda_m}(y)\quad\text{and}\quad
J_{\Lambda_{m+q_m}^c}(z)=J_{\Lambda_{m+q_m}^c}(x)\,.
$$
\end{defn}

Condition \ref{defnperiodicdecoupling} with no translation, that is $\ell\equiv 0$, is the condition
$(D)$ of Ruelle.
For the dimension $d=1$ there is a variant of  the decoupling condition.

\begin{defn}\label{defnshiftdecoupling}
A shift space $X\subset \tA^\Z$ satisfies the \emph{$1$-decoupling condition} if \\
1) there exists a function $\bar{q}\colon \N\to\N$, written $\bar{q}_m$,
such that $\lim_{m\to\infty} \bar{q}_m/m=0$;\\
2) for $m\in \N$ and  $x^-,y, x^+\in X$ there
exist $\ell^-$ and $\ell^+$, $|\ell^+|,|\ell^-|\leq \bar{q}_m$, and  $z\in X$ such that
\begin{eqnarray*}
J_{[-m,m]}(z)&=&J_{[-m,m]}(y)\\
J_{(-\infty, -m-\bar{q}_m)}(z)=
J_{(-\infty, -m-\bar{q}_m)}(T^{-\ell^-}x^-)&\text{and}&
J_{(m+\bar{q}_m,\infty)}(z)=
J_{(m+\bar{q}_m,\infty)}(T^{-\ell^+}x^+)\,.
\end{eqnarray*}
\end{defn}

If $d=1$ and $x^+=x^-=x$, then condition \ref{defnperiodicdecoupling}  implies condition \ref{defnshiftdecoupling} with
$\ell^-=\ell^+=\ell$ and $\bar{q}_m=2q_m$.
Without restricting the generality, we assume from now on that $q_m$ and $\bar{q}_m$ are monotone non-decreasing.

\begin{thm}\label{thm2}
Let $X\subset \tA^{\Z^d}$ be a shift-space satisfying the decoupling condition
\ref{defnperiodicdecoupling} and $\Phi\in\B$.   If
$\nu\in M_1(X,T)$ is a tangent functional to the pressure at $\varphi_\Phi$, then  $\nu$
is a weak Gibbs measure for $\psi=\varphi_\Phi-P(\varphi_\Phi)$.
 If $d=1$, then the same conclusion holds  under the $1$-decoupling condition
\ref{defnshiftdecoupling} and $\Phi\in\B$ such that
$$
\|\Phi\|^{\prime} :=\sum_{A\subset [-1,1]^c}\|\Phi_A\|_\infty< \infty\,.
$$
\end{thm}

\noindent
{\bf Example.\,}  Let $X$ be an irreducible sofic shift (see e.g. \cite{LM}). There exists a directed irreducible finite graph $G=(V,E)$
and a labeling of the edges $\cL\colon E\ra\tA$ so that $x\in X$ if and only if there is a bi-infinite path on $G$,
$(\ldots,e(1),e(0),e(1),\ldots)$ with $x(k)=\cL(e(k))$ for all $k\in\Z$.  Since the finite graph $G$ is irreducible,
there exists $q\in\N$ such that there is a path of length smaller than $q$ from any vertex $P\in V$ to any vertex $Q\in V$. Let $x^-,y,x^+\in X$ and $[-m,m]$ be given.
There is a  bi-infinite path $(\ldots,e^-(1),e^-(0),e^-(1),\ldots)$ which presents $x^-$,
a bi-infinite path
$(\ldots,e^\prime(1),e^\prime(0),e^\prime(1),\ldots)$ which presents $y$, and
 a  bi-infinite path $(\ldots,e^+(1),e^+(0),e^+(1),\ldots)$ which presents $x^+$.
 We have
$$
J_{[-m,m]}(y)=(y(-m),\ldots,y(m))=(\cL(e^\prime(-m)),\ldots, \cL(e^\prime(m)))\,.
$$
The path $(e^\prime(-m),\ldots, e^\prime(m))$ goes from some vertex $Q$ to some vertex $R$.
The infinite sequence $(\ldots, x^-(-(m+q)-2), x^-(-(m+q)-1))$ is presented by the infinite path
$(\ldots, e^-(-(m+q)-2),e^-(-(m+q)-1))$ which ends at some vertex $P$. Similarly the infinite sequence
$(x^+(m+q+1),x^+(m+q+2),\ldots)$ is presented by the infinite path
$(e^+(m+q+1),e^+(m+q+2),\ldots)$ which starts at some vertex $S$.
There is a path of length $q^-\leq q$ from $P$ to $Q$ and a path of length $q^+\leq q$ from $R$ to $S$. The concatenation of these  paths define a bi-infinite path presenting some $z^\prime\in X$. Then
$z:=T^{-(q-q^-)}z^\prime$ verifies the properties of definition \ref{defnshiftdecoupling} with
$\bar{q}_m=q$, $\ell^-=q-q^-$ and $\ell^+=q^+-q$.

If the shift is aperiodic, then there exists $q$ such that there exists  a path of length $q$ from $P$ to $Q$, for any $P,Q\in V$. Hence the shift verifies condition \ref{defnperiodicdecoupling} with $\ell\equiv 0$.
\qed

\begin{pro}\label{profactor}
Suppose that $X$ is a subshift that verifies  conditions
\ref{defnperiodicdecoupling} or  \ref{defnshiftdecoupling}
and that $Y$ is a factor of $X$. Then $Y$ verifies  conditions \ref{defnperiodicdecoupling} or  \ref{defnshiftdecoupling}.
\end{pro}

\prf
There is a continuous surjective map $\phi\colon X\ra Y$ such that $\phi\circ\sigma_X=\sigma_Y\circ\phi$, where $\sigma_X$ and
$\sigma_Y$ are  the shift-maps on $X$ and $Y$.
Hence $\phi$ is a sliding block-code (\cite{LM} theorem 6.2.9), i.e. there exists $\Lambda_k$ such that the $\ell^{th}$ coordinate of
$x^\prime=\phi(x)$, $x^\prime(\ell)$, depends only  on the restriction of $x$ to $\Lambda_k+\ell$.
Assume that $X$ verifies  condition \ref{defnperiodicdecoupling} with $q_m$ and
let $\Lambda_m$, $x^\prime$ and $y^\prime$ be given. Choose $x,y\in X$, such that $\phi(x)=x^\prime$ and
$\phi(y)=y^\prime$.
For $\Lambda_{m+k}$, $x$ and $y$, condition \ref{defnperiodicdecoupling} implies the existence of
$z\in X$ and $\ell\in\Lambda_{q_{m+k}}$ such that
$$
J_{\Lambda_{m+k}}(T^{-\ell}z)=J_{\Lambda_{m+k}}(y)
\quad\text{and}\quad
J_{\Lambda_{(m+k)+q_{m+k}}^c}(z)=J_{\Lambda_{(m+k)+q_{m+k}}^c}(x)\,.
$$
Let $z^\prime:=\phi(z)$. Therefore
$$
(T^{-\ell}z^\prime)(j)=(\phi(T^{-\ell}z))(j)=(\phi(y))(j)=y^\prime(j)\quad\text{for all $j\in\Lambda_m$}
$$
and
$$
z^\prime(j)=(\phi(z))(j)=(\phi(x))(j)=x^\prime(j)\quad\text{for all $j\in \Lambda_{m+2k+q_{m+k}}^c$}\,.
$$
Hence $Y$ verifies condition \ref{defnperiodicdecoupling} with $q^\prime_{m}=2k+q_{m+k}$ and $\ell^\prime=\ell$.

Assume that $X$ verifies  condition \ref{defnshiftdecoupling} with $\bar{q}_m$
and let $\Lambda_m$, ${x^\prime}^-$, $y^\prime$ and ${x^\prime}^+$ be given. Choose $x^-,y,x^+\in X$, such that $\phi(x^-)={x^\prime}^-$,
$\phi(y)=y^\prime$ and $\phi(x^+)={x^\prime}^+$. Applying condition \ref{defnshiftdecoupling} to $[-(m+k),(m+k)]$, $x^-$, $y$ and $x^+$,
one verifies as above that $Y$ satisfies condition \ref{defnshiftdecoupling}
with
$\bar{q}^\prime_m=\bar{q}_{m+k}+2k$ and
$(\ell^\prime)^\pm=\ell^\pm$.
\qed

\section{Proof  theorem \ref{thm2}}\label{section2}
\setcounter{equation}{0}

We first state two general lemmas.

\iffalse
Lemma \ref{lemq_m} allows to assume that $q_m$ is monotone increasing.

\begin{lem}\label{lemq_m}
Let $q\colon \N\to\N$ be such that $\lim_{n\to\infty} q_m/m=0$. Then $r_m:=\max\{q_j\colon j\leq m\}$ also verifies
$\lim_{n\to\infty} r_m/m=0$.
\end{lem}

\prf
Let $r_m:=q_{j_m}$. The result is true if  the sequence $\{q_m\}$ is bounded. If $\{q_m\}$ is unbounded, then
$\{r_m\}$ is a monotone diverging sequence and the same is true for $\{j_m\}$.
By hypothesis, given $\varepsilon>0$, there exists $n_\varepsilon$ such that $q_k/k\leq\varepsilon$ for all
$k\geq n_\varepsilon$. Therefore
$$
\frac{r_m}{m}=\frac{q_{j_m}}{j_m}\frac{j_m}{m}\leq \frac{q_{j_m}}{j_m}\leq\varepsilon\quad
\text{if $j_m\geq n_\varepsilon$.}
$$
\qed
\fi

\begin{lem}\label{lemestimates}
Let $\Phi\in\B$. Then
$$
\lim_{n\ra\infty}\frac{1}{|\Lambda_n|}\|W_{\Lambda_n,\Lambda_n^c}\|_\infty=0\,.
$$
\end{lem}

\prf
Let $\partial_A\Lambda:=\{j\in\Lambda\colon A+j\not\subset\Lambda\}$.
Since $\|\Phi_A\|_\infty=\|\Phi_{A+j}\|_\infty$,
$$
\|W_{\Lambda_n,\Lambda_n^c}\|_\infty\leq \sum_{\substack{A\cap\Lambda_n\not=
\emptyset\\A\cap \Lambda_n^c\not=\emptyset}} \|\Phi_A\|_\infty\leq
  \sum_{j\in\Lambda_n}\sum_{\substack{A\colon A\ni j\\A\not\subset \Lambda_n}} \|\Phi_A\|_\infty
  = \sum_{A\ni 0}|\partial_A\Lambda_n|\cdot\|\Phi_A\|_\infty\,.
$$
By the dominated convergence theorem
$$
\lim_n\sum_{A\ni 0}\frac{|\partial_A\Lambda_n|}{|\Lambda_n|}\|\Phi_A\|_\infty=
\sum_{A\ni 0}\Big(\lim_n\frac{|\partial_A\Lambda_n|}{|\Lambda_n|}\Big)\|\Phi_A\|_\infty=0\,.
$$
\qed

\begin{lem}\label{2.1.1}
Let  $n>m+q_m$, $q_m>0$, and $\ell,\ell^\prime$ such that $\Lambda_m+\ell\subset
\Lambda_{m+q_m}$ and $\Lambda_m+\ell^\prime\subset \Lambda_{m+q_m}$.
Given $\varepsilon>0$, there exists $m_\varepsilon$ so that for
$m\geq m_\varepsilon$
and $x,y\in X_n$ such that
$J_{\Lambda_n\setminus(\Lambda_{m+q_m})}(x)=
J_{\Lambda_n\setminus(\Lambda_{m+q_m})}(y)$,
$$
\frac{\exp\big(U_{\Lambda_m+\ell}(x)-U_{\Lambda_m+\ell^\prime}(y)\big)}{C_m(\varepsilon)}
\leq
\exp\big(U_{\Lambda_n}(x)-U_{\Lambda_n}(y)\big)
\leq
C_m(\varepsilon)\exp\big(U_{\Lambda_m+\ell}(x)-U_{\Lambda_m+\ell^\prime}(y)\big)
$$
where
$$
C_m(\varepsilon)=\exp|\Lambda_m|\Big(\frac{2|\Lambda_{m+q_m}\setminus\Lambda_m|}
{|\Lambda_m|}\,\|\Phi\|+2\varepsilon\Big)\,.
$$
The same results are true for $\Lambda_n\supset (\Lambda_{m+q_m}+j)$,
$(\Lambda_m+j+\ell)\subset
(\Lambda_{m+q_m}+j)$ and $(\Lambda_m+j+\ell^\prime)\subset (\Lambda_{m+q_m}+j)$.
\end{lem}

\prf
By lemma \ref{lemestimates} there exists $m_\varepsilon$ such that for all $m\geq m_\varepsilon$,
$$
 \|W_{\Lambda_m+j,(\Lambda_m+j)^c}\|_\infty  \leq\varepsilon|\Lambda_m|\,.
$$
$$
U_{\Lambda_n}(x)=U_{\Lambda_n\setminus(\Lambda_{m+q_m})}(x)+U_{\Lambda_m+\ell}(x)+
\sum_{\substack{A\subset\Lambda_n\backslash(\Lambda_m+\ell)\\A\cap(\Lambda_{m+q_m})\not=\emptyset}}
\Phi_A(x)+W_{\Lambda_m+\ell,\Lambda_n\backslash(\Lambda_m+\ell)}(x)\,.
$$
By hypothesis
$U_{\Lambda_n\setminus(\Lambda_{m+q_m})}(x)=U_{\Lambda_n\setminus(\Lambda_{m+q_m})}(y)$.
Therefore
$$
\big|\big(U_{\Lambda_n}(x)-U_{\Lambda_n}(y)\big)-
\big(U_{\Lambda_m+\ell}(x)-U_{\Lambda_m+\ell^\prime}(y)\big)\big|\leq
2|\Lambda_{m+q_m}\setminus\Lambda_m|\,\|\Phi\|+2\varepsilon\,|\Lambda_m|\,.
$$
\qed

\subsection{Proof of theorem \ref{thm2} under the decoupling condition \ref{defnperiodicdecoupling}}\label{subsectionthm1}

We  prove theorem \ref{thm2} when $\partial P(\varphi_\Phi)=\{\nu\}$ and condition \ref{defnperiodicdecoupling} holds. We define a potential $\Psi$. Let $\bar{u}\in J_m(X)$;
$$
\Psi_A(x):=
\begin{cases}
1 & \text{if $A=\Lambda_m+j$ and $J_{\Lambda_m+j}(x)=\bar u$}\\
0 & \text{otherwise.}
\end{cases}
$$
Since $\nu\in M_1(X,T)$ is a tangent functional to the convex
function $P$ at $\varphi_\Phi$,
\be\nonumber
   \frac{P(\varphi_\Phi)-P(\varphi_\Phi-t\varphi_\Psi)}{t}\leq
   \brakt{\nu}{\Psi_{\Lambda_m}}\leq
   \frac{P(\varphi_\Phi+t\varphi_\Psi)-P(\varphi_\Phi)}{t}.
\ee
Since $\nu$ is the unique tangent functional to $P$ at $\varphi$, $t\mapsto P(\varphi+t\varphi_\Psi)$ is differentiable at $t=0$ and one may interchange $\frac {d}{dt}$ and $\lim_n$  (see theorem 25.7 \cite{Ro}),
$$
 \brakt{\nu}{\Psi_{\Lambda_m}}=\left .\frac {d}{dt} P(\varphi+t\varphi_\Psi)\right|_{t=0}
=\lim_{n\ra\infty}\left .\frac {d}{dt} P_n(\Phi+t\Psi)\right|_{t=0}
$$
and
\begin{eqnarray}\label{diffsum}
\left .\frac {d}{dt} P_n(\Phi+t\Psi)\right|_{t=0}
 &=&
 \frac{1}{|\Lambda_n|}
 \frac{ \sum_{A\subset\Lambda_n}\sum_{x\in X_n}\Psi_A(x) \exp U_{\Lambda_n}(x)}
{ \sum_{x\in X_n}\exp U_{\Lambda_n}(x)}\\
&=&
 \frac{1}{|\Lambda_n|}
 \frac{\sum_{j\in\Lambda_{n-m}} \sum_{x\in X_n}\Psi_{\Lambda_m+j}(x)
 \exp U_{\Lambda_n}(x)}
{ \sum_{x\in X_n}\exp U_{\Lambda_n}(x)}\,.\nonumber
\end{eqnarray}
The key step of the proof is to obtain upper and lower bounds  independent of  $n$ and
$j\in \Lambda_{n-(m+q_m)}$ for
$$
\frac{ \sum_{x\in X_n}\Psi_{\Lambda_m+j}(x)
 \exp U_{\Lambda_n}(x)}
{ \sum_{x\in X_n}\exp U_{\Lambda_n}(x)}\,.
$$
(The terms with $j\in \Lambda_{n-m}\backslash\Lambda_{n-(m+q_m)}$ do not affect the limit $n\ra\infty$
in \eqref{diffsum}.)
For $j\in \Lambda_{n-(m+q_m)}$ and $\ell\in\Lambda_{q_m}$, let
$$
E^j_\ell(v):=\{x\in X_n\colon J_{(\Lambda_m+j)+\ell}(x)=v\}\quad\text{and}\quad
Z^j_{n,\ell}(v):=\sum_{x\in E^j_{\ell}(v)}\exp U_{\Lambda_n}(x)\,.
$$
Hence
$$
\frac{ \sum_{x\in X_n}\Psi_{\Lambda_m+j}(x) \exp U_{\Lambda_n}(x)}
{\sum_{x\in X_n}\exp U_{\Lambda_n}(x)}=\frac{Z^j_{n,0}(\bar u)}
{\sum_{x\in X_n}\exp U_{\Lambda_n}(x)}\,.
$$
For any  $\ell\in\Lambda_{q_m}$ and $j\in \Lambda_{n-(m+q_m)}$,
\begin{equation}\label{eqell}
\sum_{x\in X_n}\exp U_{\Lambda_n}(x)
=\sum_{v\in J_{(\Lambda_m+j)+\ell}(X)}Z^j_{n,\ell}(v)\,.
\end{equation}

\begin{lem}\label{2.1.2}
Let $\varepsilon>0$, $n>m+q_m$, $j\in \Lambda_{n-(m+q_m)}$  and   $\bar u\in X_{\Lambda_m+j}$.
Let
$$
K_m(\varepsilon)=|\Lambda_{q_m}||\tA|^{|\Lambda_{m+q_m}\setminus\Lambda_m|}C_m(\varepsilon)\,.
$$
Then there exists
$m_\varepsilon$ such that for all $m\geq m_\varepsilon$,
$$
\frac{Z^j_{n,0}(\bar u)}{\sum_{x\in X_n}\exp U_{\Lambda_n}(x)}\leq
\frac{|\Lambda_{q_m}|K_m(\varepsilon)\exp U_{\Lambda_m}(\bar u)}{\exp|\Lambda_m|P_m(\Phi)}
$$
and
$$
\frac{\sum_{\ell\in\Lambda_{q_m}}Z^j_{n,\ell}(\bar u)}{\sum_{x\in X_n}\exp U_{\Lambda_n}(x)}
\geq \frac{\exp U_{\Lambda_m}(\bar u)}{K_m(\varepsilon)\exp|\Lambda_m|P_m(\Phi)}\,.
$$
\end{lem}

\prf
To simplify the notations we consider the case $j=0$ since the other cases are treated exactly in the same manner. We write $Z^0_{n,\ell}=Z_{n,\ell}$ and $E^0_\ell(v)=E_\ell(v)$.  Let
$$
X_{n,m}:=J_{\Lambda_n\setminus\Lambda_{m+q_m}}(X)\,.
$$
For $\ell\in\Lambda_{q_m}$ we decompose $E_\ell(v)$ into the disjoint union of subsets $E_\ell(v,w)$, $w\in X_{n,m}$,
$$
E_\ell(v,w):=\{x\in X_n\colon J_{\Lambda_m+\ell}(x)=v\,,\; J_{\Lambda_n\setminus\Lambda_{m+q_m}}(x)=w\}\,.
$$
Hence
$$
Z_{n,\ell}(v)=\sum_{w\in X_{n,m}}\sum_{x\in E_\ell(v,w)}\exp U_{\Lambda_n}(x)\,.
$$
Notice that it may happen that $E_\ell(v,w)=\emptyset$ for some $(v,w)$. But, by the decoupling condition
\ref{defnperiodicdecoupling}, for any $(v,w)$ there exists $\ell\in\Lambda_{q_m}$ such that
$E_\ell(v,w)\not=\emptyset$.

By a slight abuse of notations, if $J_{\Lambda_m}(T^{\ell}z)=v$ we also write $J_{\Lambda_m+\ell}(z)=v$.
With this convention we define for $v\in X_m$ an auxiliary partition function $\bar{Z}_n(v)$. Let
$$
\bar{E}(v,w):=\bigcup_{\ell\in\Lambda_{q_m}}E_{\ell}(v,w)
$$
and
$$
\bar{Z}_n(v):=\sum_{w\in X_{n,m}}\sum_{x\in \bar{E}(v,w)}\exp U_{\Lambda_n}(x)\,.
$$
For any $\ell$
\begin{equation}\label{eqineqell1}
Z_{n,\ell}(v)=\sum_{w\in X_{n,m}}\sum_{x\in E_{\ell}(v,w)} \exp U_{\Lambda_n}(x)
\leq \bar{Z}_n(v)\leq\sum_{\ell^\prime\in\Lambda_{q_m}}Z_{n,\ell^\prime}(v)\,.
\end{equation}
Summing over $v\in X_m$, we get using \eqref{eqell},
\begin{equation}\label{eqineqell2}
\sum_{v\in X_m}Z_{n,\ell}(v)=
\sum_{x\in X_n}\exp U_{\Lambda_n}(x)\leq
\sum_{v\in X_m}\bar{Z}_n(v)\leq|\Lambda_{q_m}|\sum_{x\in X_n}\exp U_{\Lambda_n}(x)\,.
\end{equation}

We now compare the auxiliary partition functions $\bar{Z}_n(v)$ and $\bar{Z}_n(\bar{u})$.
Let $x\in \bar{E}(v,w)$ and $y\in \bar{E}(\bar u,w)$.
There exist $\ell$ and $\ell^\prime$ such that $x\in E_\ell(v,w)$ and $y\in E_{\ell^\prime}(\bar u,w)$.
By lemma \ref{2.1.1}
\begin{equation}\label{ub}
\exp U_{\Lambda_n}(x)
\leq
C_m(\varepsilon)
\exp(U_{\Lambda_m+\ell}(v)-U_{\Lambda_m+\ell^\prime}(\bar u))
\exp U_{\Lambda_n}(y)
\end{equation}
and
\begin{equation}\label{lb}
C_m(\varepsilon)^{-1}\exp(U_{\Lambda_m+\ell}(v)-U_{\Lambda_m+\ell^\prime}(\bar u))
\exp U_{\Lambda_n}(y)\leq \exp U_{\Lambda_n}(x)\,.
\end{equation}
By translation invariance of $\Phi$,  if $J_{\Lambda_m+\ell}(x)=v$ and $J_{\Lambda_m}(x^\prime)=v$, then
$$
U_{\Lambda_m+\ell}(v)\equiv U_{\Lambda_m+\ell}(x)=U_{\Lambda_m}(x^\prime)\equiv U_{\Lambda_m}(v)\,.
$$
Similarly,
$U_{\Lambda_m+\ell^\prime}(\bar u)=U_{\Lambda_m}(\bar u)$.
For any $v$, the cardinality of $E_\ell(v,w)$ is smaller than $|\tA|^{|\Lambda_{m+q_m}\setminus\Lambda_m|}$ and that of
$\bar{E}(v,w)$ smaller than $|\Lambda_{q_m}||\tA|^{|\Lambda_{m+q_m}\setminus\Lambda_m|}$.
Summing \eqref{ub} over $y$ and then over $x$,
$$
\sum_{x\in \bar{E}(v,w)}\exp U_{\Lambda_n}(x)\leq
K_m(\varepsilon)\exp( U_{\Lambda_m}(v)-U_{\Lambda_m}(\bar u))
\sum_{y\in \bar{E}(\bar u,w)}\exp U_{\Lambda_n}(y)\,.
$$
Summing \eqref{lb} over $x$ and then over $y$,
$$
\exp(U_{\Lambda_m}(v)-U_{\Lambda_m}(\bar u))
\frac{\sum_{y\in \bar{E}(u,w)}\exp U_{\Lambda_n}(y)}{K_m(\varepsilon)}
\leq
\sum_{x\in \bar{E}(v,w)}\exp U_{\Lambda_n}(x)\,.
$$
Summing over $w\in X_{n,m}$ we get
$$
{\rm e}^{U_{\Lambda_m}(v)}{\rm e}^{-U_{\Lambda_m}(\bar u)}
\frac{\bar{Z}_n(\bar u)}{K_m(\varepsilon)}
\leq
\bar{Z}_n(v)
\leq
{\rm e}^{U_{\Lambda_m}(v)}{\rm e}^{-U_{\Lambda_m}(\bar u)}
K_m(\varepsilon)\bar{Z}_n(\bar u)\,.
$$
Finally we sum over $v$ and use \eqref{eqineqell1} and \eqref{eqineqell2}.
We get by  \eqref{eqineqell2} and then \eqref{eqineqell1}
\begin{eqnarray*}
\sum_{v}\frac{Z_{n,0}(v)}{\sum_{\ell\in\Lambda_{q_m}}Z_{n,\ell}(\bar u)}
&=&
\frac{\sum_{x\in X_n}\exp U_{\Lambda_n}(x)}{\sum_{\ell\in\Lambda_{q_m}}Z_{n,\ell}(\bar u)}\\
&\leq &
\frac{\sum_v\bar{Z}_n(v)}{\sum_{\ell\in\Lambda_{q_m}}Z_{n,\ell}(\bar u)}\leq
\frac{\sum_v\bar{Z}_n(v)}{\bar{Z}_n(\bar u)}\leq
K_m(\varepsilon)\frac{\sum_v\exp U_{\Lambda_m}(v)}{\exp U_{\Lambda_m}(\bar u)}\,,
\end{eqnarray*}
so that
$$
\frac{\sum_{\ell\in\Lambda_{q_m}}Z_{n,\ell}(\bar u)}{\sum_{x\in X_n}\exp U_{\Lambda_n}(x)}
\geq
\frac{\exp U_{\Lambda_m}(\bar u)}{K_m(\varepsilon)\exp|\Lambda_m|P_m(\Phi)}\,.
$$
Similarly,  $$
\sum_{v}\frac{\sum_{\ell\in\Lambda_{q_m}}Z_{n,\ell}(v)}{Z_{n,0}(\bar u)}
=
\frac{|\Lambda_{q_m}|\sum_{x\in X_n}\exp U_{\Lambda_n}(x)}{Z_{n,0}(\bar u)}
\geq
\frac{\sum_v\bar{Z}_n(v)}{\bar{Z}_n(\bar u)}\geq
\frac{\sum_v\exp U_{\Lambda_m}(v)}{K_m(\varepsilon)\exp U_{\Lambda_m}(\bar u)}\,,
$$
so that
$$
\frac{Z_{n,0}(\bar u)}{\sum_{x\in X_n}\exp U_{\Lambda_n}(x)}
\leq
\frac{|\Lambda_{q_m}|K_m(\varepsilon)\exp U_{\Lambda_m}(\bar u)}{\exp|\Lambda_m|P_m(\Phi)}\,.
$$
\qed

\begin{lem}\label{2.1.3}
Uniformly in $z\in X$,
$$
\lim_{m\ra\infty}\Big(\brakt{\cE_m(z)}{\varphi_\Phi}-\frac{U_{\Lambda_m}(z)}{|\Lambda_m|}\Big)=0\,.
$$
\end{lem}
\prf
\begin{eqnarray*}
\Big| |\Lambda_m|\brakt{\cE_m(z)}{\varphi_\Phi}-U_{\Lambda_m}(z)\Big|
&=&
\Big|\sum_{k\in\Lambda_m}\varphi_\Phi(T^kz)-\sum_{A\subset \Lambda_m}\Phi_A(z)\Big|\\
&=&
\Big|\sum_{k\in\Lambda_m}\Big(\sum_{A\ni k}\frac{\Phi_A(z)}{|A|}-
\sum_{\substack{A\ni k\\A\subset\Lambda_m}}\frac{\Phi_A(z)}{|A|}\Big)\Big|\\
&\leq&
\Big|\sum_{k\in\Lambda_m}\sum_{\substack{A\ni k\\A\cap\Lambda^c_m\not=\emptyset}}
\frac{\Phi_A(z)}{|A|}\Big|\\
&\leq&
\Big|\sum_{\substack{A\cap\Lambda_m\not=\emptyset\\A\cap\Lambda^c_m\not=\emptyset}}
\frac{|A\cap\Lambda_m|\Phi_A(z)}{|A|}\Big|\leq \|W_{\Lambda_m,\Lambda_m^c}\|_\infty\,.
\end{eqnarray*}
\qed
\noindent
Let $m\geq m_\varepsilon$ and $z\in X$, $J_m(z)=\bar u$.
Taking the limit $n\ra\infty$ in \eqref{diffsum}, one gets from lemma \ref{2.1.2}
and
$$
\sum_{j\in\Lambda_{n-m}}Z_{n,0}^j(\bar u)\leq \sum_{j\in\Lambda_{n-m-q_m}}
\sum_{\ell\in \Lambda_{q_m}}Z^j_{n,\ell}(\bar u)
\leq
|\Lambda_{q_m}|\sum_{j\in\Lambda_{n-m}}Z^j_{n,0}(\bar u)\,,
$$
$$
\frac{\exp U_{\Lambda_m}(z)}{|\Lambda_{q_m}|K_m(\varepsilon)\exp|\Lambda_m|P_m(\Phi)}
\leq
 \brakt{\nu}{\Psi_{\Lambda_m}}
 \iffalse
\frac{\sum_{x\in X_n}\Psi_{\Lambda_m}(x)\exp U_{\Lambda_n}(x)}
{ \sum_{y\in X_n}\exp U_{\Lambda_n}(y)}
\fi
\leq
\frac{|\Lambda_{q_m}|K_m(\varepsilon)\exp U_{\Lambda_m}(z)}{\exp|\Lambda_m|P_m(\Phi)}\,.
$$
Therefore, by lemma \ref{2.1.3}, given $\delta>0$ there exists $N_\delta$ such that
\begin{equation}\label{final}
\sup_{z\in X}\Big|\frac{1}{|\Lambda(m)|}\ln\nu(B_{m}(z))-\brakt{\cE_m(z)}{(\varphi_\Phi-P(\varphi_\Phi)}\Big|\leq \delta\quad\forall m\geq N_\delta\,.
\end{equation}
This proves the theorem when $\partial P(\varphi_\Phi)=\{\nu\}$.

The pressure is convex and continuous on the Banach space $\B$.
To prove the theorem in the general case we apply a theorem of Mazur and a theorem of Lanford and Robinson
(see  \cite{Ru} appendix A.3.7).
The set of $\Phi$ such that $\partial P(\varphi_\Phi)=\{\nu\}$ has a unique element $\nu$ is residual (theorem of Mazur). Let $\nu\in\partial P(\varphi_\Phi)$ be such that there exist
$\Phi_k\in\B$ with $\partial P(\varphi_{\Phi_k})=\{\nu_k\}$, $\lim_k\Phi_k=\Phi$
and $\lim_k\nu_k=\nu$.
One can choose $N_\delta$ so that \eqref{final} holds for $\nu_k$ and $\varphi_{\Phi_k}$, uniformly in $k$.
Hence \eqref{final} holds also for such $\nu\in\partial P(\varphi_\Phi)$ and $\varphi_\Phi$.
It also holds for
$\mu$ in the convex hull of such $\nu$ and for all
$\rho\in \partial P(\varphi)$ which are limits of such $\mu$. Hence \eqref{final} is true on the
weak-closed convex hull of such $\nu$, which coincides with
$\partial P(\varphi_\Phi)$ (theorem of Lanford and Robinson).
\qed

\subsection{Proof of theorem \ref{thm2} under the $1$-decoupling condition \ref{defnshiftdecoupling}}\label{subsectionthm2}

The scheme of the proof is the same as before.
Let $n> m+2\bar{q}_m$  and $|j|< n-(m+3\bar{q}_m)$ be given. To simplify the notation, since below $n$, $m$ and $j$ are fixed,
we write the complement of
$(\Lambda_{m+2\bar{q}_m}+j)$ in $\Lambda_{n-\bar{q}_m}$ as
$$
\Lambda^-\cup\Lambda^+:=\Lambda_{n-\bar{q}_m}\setminus(\Lambda_{m+2\bar{q}_m}+j)\quad \text{where}\;
\Lambda^-= [-n+\bar{q}_m,-(m+2\bar{q}_m)+j)\,.
$$
By our choice of $j$, $\Lambda^-\not=\emptyset$ and $\Lambda^+\not=\emptyset$.  We set
$$
X^{-}:=J_{\Lambda^-}(X)\,,\; X^{+}:=J_{\Lambda^+}(X)\,,\; E(v):=\{x\in X_n\colon J_{\Lambda_m+j}(x)=v\}\,,
$$
and
$$
Z_n(v):=\sum_{x\in E(v)}\exp U_{\Lambda_n}(x)\,.
$$
For $\Lambda\subset\Z$, we also denote by $T^\ell$ the map from $J_\Lambda(X)$ to $J_{\Lambda-\ell}(X)$ defined by
$$
w\mapsto (T^\ell w)(k):=w(k+\ell)\,,\; k\in(\Lambda-\ell)\,.
$$
Let $|\ell^-|,|\ell^+|\leq\bar{q}_m$ be given;  we decompose the set $E(v)$ as the disjoint union of  sets
$$
E_{\ell^-,\ell^+}(v,w^-,w^+):=
\{x\in E(v)\colon
J_{\Lambda^-+\ell^-}(x)=T^{-\ell^-}w^-\,,\; J_{\Lambda^++\ell^+}(x)=T^{-\ell^+}w^+\}\,,
$$
where $w^-\in X^-$ and $w^+\in X^+$.
Hence for fixed $\ell^-,\ell^+$
$$
Z_n(v)
=
\sum_{w^-,w^+}\sum_{x\in E_{\ell^-,\ell^+}(v,w^-,w^+)}
\exp U_{\Lambda_n}(x)\,.
$$
The cardinality of each  $E_{\ell^-,\ell^+}(v,w^-,w^+)$ is at most $|\tA|^{6\bar{q}_m}$. The set
$E_{\ell^-,\ell^+}(v,w^-,w^+)$ may be empty, but condition
\ref{defnshiftdecoupling} implies that $E_{\ell^-,\ell^+}(v,w^-,w^+)\not=\emptyset$ for at least one pair $(\ell^-,\ell^+)$ for each
triple $(v,w^-,w^+)$. Indeed, for any $w^-\in X^-$ and $w^+\in X^+$ there exist $x^-\in X$ and $x^+\in X$ such that
$J_{\Lambda^-}(x^-)=w^-$ and $J_{\Lambda^+}(x^+)=w^+$. Condition \ref{defnshiftdecoupling} implies the existence of $z\in X$ such that  $J_{[-m,m]}(z)=v$ and
$$
J_{(-\infty, -m-\bar{q}_m)}(z)=
J_{(-\infty, -m-\bar{q}_m)}(T^{-\ell^-}x^-)\quad\text{and}\quad
J_{(m+\bar{q}_m,\infty)}(z)=
J_{(m+\bar{q}_m,\infty)}(T^{-\ell^+}x^+)\,.
$$
Since ${\Lambda^-+\ell^-}\subset [-n, -m-\bar{q}_m)$ and $\Lambda^++\ell^+\subset (m+\bar{q}_m,n]$, we have
$$
J_{\Lambda^-+\ell^-}(z)=J_{\Lambda^-+\ell^-}(T^{-\ell^-}x^-)=T^{-\ell^-}w^-\,,\;
J_{\Lambda^++\ell^+}(z)=J_{\Lambda^++\ell^+}(T^{-\ell^+}x^+)=T^{-\ell^+}w^+\,.
$$
We introduce an auxiliary partition function  $\bar{Z}_n(v)$.  Let
$$
\bar{E}(v,w^-,w^+):=\bigcup_{\ell^-,\ell^+} E_{\ell^-,\ell^+}(v,w^-,w^+)\quad\text{and}\quad
\bar{Z}_n(v):= \sum_{w^-,w^+}\sum_{x\in \bar{E}(v,w^-,w^+)}\exp U_{\Lambda_n}(x)\,.
$$
We have
\begin{equation}\label{cardinality}
1\leq |\bar{E}(v,w^-,w^+)|\leq (2\bar{q}_m+1)^2|\tA|^{6\bar{q}_m}
\end{equation}
and
\begin{equation}\label{trick}
Z_n(v)\leq \bar{Z}_n(v)\leq (2\bar{q}_m+1)^2Z_n(v)\,.
\end{equation}
The next step is to compare $\bar{Z}_n(v)$ and $\bar{Z}_n(\bar u)$. Let $[-n,n]$ be the disjoint union of $\Lambda^1$, $\Lambda^2$ and
$\Lambda^3$. Using this decomposition of $[-n,n]$, the energy $U_{[-n,n]}$ is equal to
\begin{equation}\label{eqdecomp}
\underbrace{\sum_{A\subset\Lambda^1}\Phi_A}_{\text{interaction inside $\Lambda^1$}}
+
\underbrace{\sum_{\substack{A\subset\Lambda^1\cup\Lambda^2\\A\cap\Lambda^1\not=\emptyset\,
A\cap(\Lambda^2)\not=\emptyset}}\Phi_A}_{\substack{\text{interaction between }\\
\text{$\Lambda^1$ and $\Lambda^2$}}}
+
\underbrace{\sum_{A\subset \Lambda^2}\Phi_A}_{\text{interaction inside $\Lambda^2$}}
+
\underbrace{\sum_{\substack{A\subset [-n,n]\\A\cap\Lambda^3\not=\emptyset}}\Phi_A}_{
\substack{\text{interaction inside $\Lambda^3$ and }\\\text{between $\Lambda^3$ and $[-n,n]\setminus{\Lambda^3}$ }}}\,.
\end{equation}
The second sum in \eqref{eqdecomp} is $W_{\Lambda^1,\Lambda^2}$ and the fourth sum is bounded by
$$
\big\|\sum_{\substack{A\subset [-n,n]\\A\cap\Lambda^3\not=\emptyset}}\Phi_A\big\|\leq
\sum_{k\in\Lambda^3}\sum_{A\ni k}\|\Phi_A\|= |\Lambda^3|\|\Phi\|\,.
$$
Let $x \in \bar{E}(v,w^-,w^+)$,  $y\in \bar{E}(\bar{u},w^-,w^+)$; hence
$x\in   E_{\ell^-,\ell^+}(v,w^-,w^+)$ for some $(\ell^-,\ell^+)$, and
$y\in   E_{{\ell^\prime}^-,{\ell^\prime}^+}(\bar{u},w^-,w^+)$ for some
$({\ell^\prime}^-,{\ell^\prime}^+)$.
We write $U_{[-n,n]}(x)$ using the identity \eqref{eqdecomp}
 with
 $$
 \Lambda^1=[-m,m]+j\quad \Lambda^2=(\Lambda^-+\ell^-)\cup(\Lambda^++\ell^+)\quad
 \Lambda^3=[-n,n]\setminus(\Lambda^1\cup\Lambda^2)\equiv \underline{\Lambda}\,,
 $$
and  $U_{[-n,n]}(y)$ using the identity \eqref{eqdecomp}, where this time
 $$
 \Lambda^1=[-m,m]+j\quad \Lambda^2=(\Lambda^-+{\ell^\prime}^-)\cup(\Lambda^++{\ell^\prime}^+)\quad
 \Lambda^3=[-n,n]\setminus(\Lambda^1\cup\Lambda^2)\equiv \underline{\Lambda}^\prime\,.
 $$
We have
$$
|\sum_{\substack{A\subset [-n,n]\\A\cap\underline{\Lambda}\not=\emptyset}}\Phi_A(x)|\leq
|\underline{\Lambda}|\|\Phi\|
\quad\text{and}\quad
|\sum_{\substack{A\subset [-n,n]\\A\cap\underline{\Lambda}^\prime\not=\emptyset}}\Phi_A(y)|\leq
|\underline{\Lambda}^\prime|\|\Phi\|\,,
$$
with $|\underline{\Lambda}^\prime|=|\underline{\Lambda}|=6\bar{q}_m$.
Since
$$
\sum_{A\subset (\Lambda^-+\ell^-)\cup(\Lambda^++\ell^+)}\Phi_A(x)=
\sum_{A\subset \Lambda^-+\ell^-}\Phi_A(x)+\sum_{A\subset \Lambda^++\ell^+}\Phi_A(x)+
W_{(\Lambda^-+\ell^-),(\Lambda^++\ell^+)}(x)\,,
$$
and by definition of $x$ and $y$
$$
\sum_{A\subset \Lambda^-+\ell^-}\Phi_A(x)=\sum_{A\subset \Lambda^-+{\ell^\prime}^-}\Phi_A(y)\quad
\text{and}\quad
\sum_{A\subset \Lambda^++\ell^+}\Phi_A(x)=\sum_{A\subset \Lambda^++{\ell^\prime}^+}\Phi_A(y)\,,
$$
we get
$$
\Big|\sum_{A\subset (\Lambda^-+\ell^-)\cup(\Lambda^++\ell^+)}\Phi_A(x)-
\sum_{A\subset (\Lambda^-+{\ell^\prime}^-)\cup(\Lambda^++{\ell^\prime}^+)}\Phi_A(y)\Big|\leq
2\|\Phi\|^\prime\,.
$$
Given $\varepsilon>0$, if $m$ is large enough,
$$
\|W_{\Lambda_m+j,(\Lambda^-+{\ell^\prime}^-)\cup(\Lambda^++{\ell^\prime}^+)}\|
\leq
\|W_{\Lambda_m+j,\Lambda_n\setminus(\Lambda_m+j)}\|\leq\varepsilon|\Lambda_m|\,.
$$
The same bound holds for $\|W_{\Lambda_m+j,(\Lambda^-+\ell_-)\cup(\Lambda^++\ell_+)}\|$.
Hence,
$$
\big|\big(U_{\Lambda_n}(x)-U_{\Lambda_m+j}(x)\big)-\big(U_{\Lambda_n}(y)-U_{\Lambda_m+j}(y)\big)\big|
\leq
2(\varepsilon|\Lambda_m|+\|\Phi\|^\prime+6\bar{q}_m\|\Phi\|)\,,
$$
and
\begin{equation}\label{est}
\frac{{\rm e}^{U_{\Lambda_m+j}(x)-U_{\Lambda_m+j}(y)}}{C^\prime_m(\varepsilon)}\leq
{\rm e}^{U_{\Lambda_n}(x)-U_{\Lambda_n}(y)}\leq  C^\prime_m(\varepsilon){\rm e}^{U_{\Lambda_m+j}(x)-U_{\Lambda_m+j}(y)}
\end{equation}
with
$$
C^\prime_m(\varepsilon)=\exp{2|\Lambda_m|\Big(\varepsilon+\frac{6\bar{q}_m|\Phi\|+\|\Phi\|^\prime}{|\Lambda_m|}}\Big)\,.
$$
Let
$$
\bar{K}_m(\varepsilon):=(2\bar{q}_m+1)^2|\tA|^{6\bar{q}_m}C^\prime_m(\varepsilon)\,.
$$
Summing \eqref{est} over $x\in \bar{E}(v,w^-,w^+)$ and $y\in \bar{E}(\bar{u},w^-,w^+)$, taking into account \eqref{cardinality}, and then summing over $(w^-,w^+)\in X^-\times X^+$, we get
$$
{\rm e}^{U_{\Lambda_m+j}(v)}{\rm e}^{-U_{\Lambda_m+j}(\bar{u})}\frac{\bar{Z}_n(\bar{u})}{\bar{K}_m(\varepsilon)}\leq
\bar{Z}_n(v)\leq \bar{K}_m(\varepsilon){\rm e}^{U_{\Lambda_m+j}(v)}{\rm e}^{-U_{\Lambda_m+j}(\bar{u})}\bar{Z}_n(\bar{u})\,.
$$
Using \eqref{trick} we get by summing over $v\in X_m$,
$$
\frac{\exp U_{\Lambda_m+j}(\bar u)}{(2\bar{q}_m+1)^2\bar{K}_m(\varepsilon)\exp |\Lambda_m|P_m(\Phi)}\leq
\frac{Z_n(\bar u)}{\sum_v Z_n(v)}\leq \frac{(2\bar{q}_m+1)^2\bar{K}_m(\varepsilon)\exp U_{\Lambda_m+j}(\bar u)}
{\exp |\Lambda_m| P_m(\Phi)}\,.
$$
The rest of the proof is identical to the proof of theorem \ref{thm2} in subsection \ref{subsectionthm1}
when there is a unique tangent functional to $P$ at $\Phi$.
The general case follows as in subsection \ref{subsectionthm1} by noticing that $P$ is convex and continuous on the Banach space
$\{\Phi\in\B\colon \vvvert \Phi\vvvert=\|\Phi\|+\|\Phi\|^\prime<\infty\}$.
\qed

\vspace*{1cm}
\noindent
{\bf Acknowledgments.\,}
C.-E. Pfister thanks N. Cuneo for correspondence.

\end{document}